\documentclass{amsart}
\usepackage{geometry,amsmath,amssymb,graphicx,bbm}
\renewcommand{\email}[1]{{\textit{Email:} \texttt{#1}}}
\newcommand{\mathd}{\mathrm{d}}
\newcommand{\matheuler}{\gamma}
\newcommand{\nin}{\not\in}
\newcommand{\tmem}[1]{{\em #1\/}}
\newcommand{\tmop}[1]{\ensuremath{\operatorname{#1}}}
\newcommand{\tmsamp}[1]{\textsf{#1}}
\newcommand{\tmtextit}[1]{{\itshape{#1}}}
\newcommand{\tmtexttt}[1]{{\ttfamily{#1}}}

\begin{document}

\title{Two New Zeta Constants: Fractal String, Continued Fraction, and
Hypergeometric Aspects of the Riemann Zeta Function}\author{Stephen
Crowley}\thanks{\email{stephen.crowley@mavs.uta.edu}}\maketitle

\begin{abstract}
  The Riemann zeta function at integer arguments can be written as an infinite
  sum of certain hypergeometric functions and more generally the same can be
  done with polylogarithms, for which several zeta functions are a special
  case. An analytic continuation formula for these hypergeometric functions
  exists and is used to derive some infinite sums which allow the zeta
  function at integer arguments $n$ to be written as a weighted infinite sum
  of hypergeometric functions at $n - 1$. The form might be considered to be a
  shift operator for the Riemann zeta function which leads to the curious
  values $\zeta_{}^F (0) = I_0 (2) - 1$ and $\zeta_{}^F (1) = \tmop{Ei} (1) -
  \gamma$ which involve a Bessel function of the first kind and an exponential
  integral respectively and differ from the values $\zeta \left( 0 \right) = -
  \frac{1}{2}$ and $\zeta \left( 1 \right) = \infty$ given by the usual method
  of continuation. Interpreting these ``hypergeometrically continued'' values
  of the zeta constants in terms of reciprocal common factor probability we
  have $\zeta^F \left( 0 \right)^{- 1} \cong 78.15\%$ and $\zeta^F \left( 1
  \right)^{- 1} \cong 75.88\%$ which contrasts with the standard known values
  for sensible cases like $\zeta \left( 2 \right)^{- 1} = \frac{6}{\pi^2}
  \cong 60.79\%$ and $\zeta \left( 3 \right)^{- 1} \cong 83.19\%$. The
  combinatorial definitions of the Stirling numbers of the second kind, and
  the $2$-restricted Stirling numbers of the second kind are recalled because
  they appear in the differential equation satisfied by the hypergeometric
  representation of the polylogarithm. The notion of fractal strings is
  related to the (chaotic) Gauss map of the unit interval which arises in the
  study of continued fractions, and another chaotic map is also introduced
  called the ``Harmonic sawtooth'' whose Mellin transform is the
  (appropritately scaled) Riemann zeta function. These maps are within the
  family of what might be called ``deterministic chaos''. Some number
  theoretic definitions are also recalled.
\end{abstract}

{\tableofcontents}

\section{The Zeta Function}

\subsection{Riemann's $\zeta (t)$ Function}

{\tmem{Riemann's zeta function}}, named after Bernhard Riemann(1826-1866), is
defined by
\begin{equation}
  \begin{array}{lll}
    \zeta (t) & = \sum_{n = 1}^{\infty} n^{- t} & \forall \{t \in \mathbbm{C}:
    \mathcal{R}(t) > 1\}\\
    & = \frac{1}{1 - 2^{- t}} \sum_{n = 0}^{\infty} (2 n + 1)^{- t} & \forall
    \{t \in \mathbbm{C}: \mathcal{R}(t) > 1\}\\
    & = \left( 1 - 2^{- t} \right) \eta \left( t \right) & \forall \{t \in
    \mathbbm{C}: \mathcal{R}(t) > 0\}
  \end{array}  \label{zeta}
\end{equation}
where \tmtexttt{{\tmsamp{}}}$\mathcal{R}(t)$ and $\mathcal{I}(t)$ denote real
and imaginary parts of $t$ respectively and $\eta \left( t \right)$ is the
{\tmem{Dirichlet eta function}}, also known as the {\tmem{alternating zeta
function}}, named after Johann Dirichlet(1805-1859)
\begin{equation}
  \begin{array}{lll}
    \eta \left( t \right) & = \sum_{n = 1}^{\infty} \frac{\left( - 1
    \right)^{n - 1}}{n^t_{}} & \forall \{t \in \mathbbm{C}: \mathcal{R}(t) >
    0\}\\
    & = \frac{1}{\Gamma \left( s \right)} \int x^{s - 1}  \frac{1}{e^x + 1}
    \mathd x & \forall \{t \in \mathbbm{C}: \mathcal{R}(t) > 0\}
  \end{array}
\end{equation}
where the integral is a Mellin transform of $\left( e^x + 1 \right)^{- 1}$.
The function $\zeta (t)$ is analytic and uniformly convergent when
$\mathfrak{R}(t) > 1$ or $\mathfrak{R}(t) > 0$ when using the eta function
form. The only singularity of $\zeta (t)$ is at $t = 1$ where it becomes the
divergent harmonic series. The reflection functional equation
{\cite[13.151]{whittaker-watson}} which relates $\zeta (t)$ to $\zeta (1 - t)$
is given by
\begin{equation}
  \zeta (t) \pi^{- t} 2^{1 - t} \Gamma (t) \cos \left( \frac{t \pi}{2} \right)
  = \zeta (1 - t) \label{reflection}
\end{equation}
The interpretation of zeta in terms of frequentist probability is that given
$n$ integers chosen at random, the probability that no common factor will
divide them all is $\zeta \left( n \right)^{- 1}$. In other words, given an
array $i$ of $n$ random intgers, $\zeta \left( n \right)^{- 1}$ is the
probabability that $\gcd \left( i_1, i_2, \ldots, i_n \right) = 1$ where
$\gcd$is the {\tmem{greatest common denominator}} function. So for example,
the probability that a pair of randomly chosen integers is coprime is $\zeta
\left( 2 \right)^{- 1} = \frac{6}{\pi^2} \cong 60.79\%$, and the probability
that a triplet of randomly chosen integers is relatively prime is $\zeta
\left( 3 \right)^{- 1} \cong 83.19\%$.
{\cite{riemann}}{\cite[13.1]{whittaker-watson}}{\cite[1.4]{edwardszeta}}

\subsubsection{The Generalized Hurwitz Zeta Function $\zeta (t, a)$}

A more general function which includes Riemann's Zeta function was defined by
A. Hurwitz.
\begin{equation}
  \begin{array}{ll}
    \zeta (t, a) & =
  \end{array} \sum_{n = 0}^{\infty} (n + a)^{- t} \label{hurwitz}
\end{equation}
Notice that the summation starts at $n = 0$ whereas Riemann's starts at $n =
1$. It is apparent that $\zeta (t)$ is a special case of $\zeta (t, a)$ where
\begin{equation}
  \zeta (t) = \sum_{n = 1}^{\infty} n^{- t} = \zeta (t, 1) = \sum_{n =
  0}^{\infty} (n + 1)^{- t} \label{hurwitzriemann}
\end{equation}
{\cite{hurwitzeta}}{\cite[13.11]{whittaker-watson}}

\subsubsection{Hypergeometric Representations of the Lerch Transcendent: $\Phi
(z, t, v)$}

The Lerch transcendent $\Phi (z, t, v)$ {\cite[1.11]{htf1}} is a further
generalization of the Hurtwitz zeta function
\begin{equation}
  \begin{array}{ll}
    \Phi (z, t, v) & = \sum_{n = 0}^{\infty} \frac{z^n}{(v + n)^t}
  \end{array}
\end{equation}
valid $\forall |z| < 1$ or $\{|z| = 1 : \mathcal{R}(t) > 1\}$ which is related
to $\zeta (t, v)$ and $\zeta (t)$ by
\begin{equation}
  \begin{array}{lll}
    \Phi (1, t, v) & = \zeta (t, v) & \\
    \Phi (1, t, 1) & = \zeta (t) & \\
    \Phi (1, t, 1 / 2) & = \zeta (t, 1 / 2) & = (2^t - 1) \zeta (t)
  \end{array} \label{phizeta}
\end{equation}
When $t = 1$ the Lerch transcendent reduces to
\begin{equation}
  \begin{array}{ll}
    \Phi (z, 1, v) & = \frac{_2 F_1 \left( \begin{array}{ll}
      1 & v\\
      1 + v & 
    \end{array} | z \right)}{v}
  \end{array}
\end{equation}
and when $n \in \{0, 1, 2, \ldots .\}$, $\Phi (z, n, v)$ has the
hypergeometric representation {\cite{hyperzeta2}}
\begin{equation}
  \begin{array}{ll}
    \Phi (z, n, v) & =
  \end{array} v^{- n} _{n + 1} F_n \left( \begin{array}{ll}
    1 & \vec{v}_n\\
    & \overrightarrow{1 + v}_n
  \end{array} | z \right)
\end{equation}
yielding
\begin{equation}
  \begin{array}{ll}
    \zeta (n, v) & = v^{- n} _{n + 1} F_n \left( \begin{array}{ll}
      1 & \vec{v}_n\\
      & \overrightarrow{1 + v}_n
    \end{array} \right)
  \end{array}
\end{equation}
and
\begin{equation}
  \begin{array}{ll}
    \zeta (n) & = \left( \frac{2^n}{2^n - 1} \right) _{n + 1} F_n \left(
    \begin{array}{ll}
      1 & \overrightarrow{1 / 2}_n\\
      & \overrightarrow{3 / 2}_n
    \end{array} \right)\\
    & =_{n + 1} F_n \left( \begin{array}{l}
      \vec{1}_{n + 1}\\
      \vec{2}_n
    \end{array} \right)
  \end{array}
\end{equation}
and thus due to (\ref{zeta}) and (\ref{phizeta}) we have the hypergeometric
transformation
\begin{equation}
  _{n + 1} F_n \left( \begin{array}{ll}
    1 & \overrightarrow{1 / 2}_n\\
    & \overrightarrow{3 / 2}_n
  \end{array} \right) = (1 - 2^{- n})_{n + 1} F_n \left( \begin{array}{l}
    \vec{1}_{n + 1}\\
    \vec{2}_n
  \end{array} \right)
\end{equation}
where the argument absent in $_p F_q$ is assumed to be 1 and the symbol
$\vec{c}_n$ denotes a parameter vector of length $n$ where each element is
equal to $c$ (e.g. $\vec{5}_3 = [5, 5, 5]$).

\subsubsection{The Hypergeometric Polylogarithm}

The polylogarithm, also known as Jonqui\`ere's function, is defined $\forall
|t| \leqslant 1, n \in \{0, 1, 2, \ldots\}$ by
\begin{equation}
  \begin{array}{lll}
    \tmop{Li}_n (t) & = \sum_{k = 1}^{\infty} \frac{t^k}{k^n} & \\
    & =_{n + 1} F_n \left( \begin{array}{l}
      \vec{1}_{n + 1}\\
      \vec{2}_n
    \end{array}   | t \right) t  & \\
    & = t \sum_{k = 0}^{\infty} \frac{t^k}{k!} \frac{\prod_{i = 1}^{n + 1}
    (1)_k}{\prod_{j = 1}^n (2)_k} & \\
    & = t \sum_{k = 0}^{\infty} t^k  \frac{(1)_k^{n + 1}}{(2)_k^n} & \\
    & = t \sum_{k = 0}^{\infty} t^k \frac{\Gamma (k + 1)^n}{\Gamma (k + 2)^n}
    & \\
    & = t \sum_{k = 0}^{\infty} \frac{t^k}{(1 + k)^n} & 
  \end{array} \label{polylog}
\end{equation}
The hypergeometric representation (\ref{hypergeom}) of $\tmop{Li}_n (t) =_{n +
1} F_n \left( \begin{array}{l}
  a_1 \ldots a_{n + 1}\\
  b_1 \ldots b_n
\end{array} | t \right) t = \tmop{Li}_n^F (t$) where $a_1 \ldots a_{n + 1} =
\vec{1}_{n + 1}$ and $b_1 \ldots b_n = \vec{2}_n$ \ is
\tmtextit{nearly-poised} \tmtextit{of the first kind} {\cite[2.1.1]{Fslater}}
since $a_1 + b_1 = 3 = \ldots = a_n + b_n = 3$. \ \ The notation
$\tmop{Li}_n^F (t)$ refers specifically the hypergeometric form of
$\tmop{Li}_n (t$). The derivatives and integrals of $\tmop{Li}_n (t$) satisfy
the recurrence relations
\begin{equation}
  \begin{array}{ll}
    \frac{\mathd}{\mathd t} \tmop{Li}_n (t) & = \frac{\tmop{Li}_{n - 1}
    (t)}{t}\\
    \frac{\mathd}{\mathd t} _{n + 1} F_n \left( \begin{array}{l}
      \vec{1}_{n + 1}\\
      \vec{2}_n
    \end{array}   | t \right) t & =_n F_{n - 1} \left( \begin{array}{l}
      \vec{1}_{n + 1}\\
      \vec{2}_n
    \end{array}   | t \right)
  \end{array}
\end{equation}
\begin{equation}
  \begin{array}{ll}
    \int_0^t \frac{\tmop{Li}_n (s)}{s} \mathd s & = \tmop{Li}_{n + 1} (t)\\
    \int_0^t {_{n + 1} F_n} \left( \begin{array}{l}
      \vec{1}_{n + 1}\\
      \vec{2}_n
    \end{array}   | s \right) \mathd s & ={_{n + 2} F_{n + 1}} \left(
    \begin{array}{l}
      \vec{1}_{n + 2}\\
      \vec{2}_{n + 1}
    \end{array}   | t \right) t
  \end{array}
\end{equation}
and the reflection functional equation for $\tmop{Li}_n (1) = \zeta (n)$ is
\begin{equation}
  \begin{array}{ll}
    \tmop{Li}_n (1) & = \frac{\tmop{Li}_n (- 1)}{(2^{1 - n} - 1)}
  \end{array} \label{lireflect}
\end{equation}
$\tmop{Li}_n^F (t)$ \ is seen to be ($n - 1$)-balanced (\ref{balance}) with
the trivial calculation
\begin{equation}
  \sum_{k = 1}^n 2 - \sum_{k = 1}^{n + 1} 1 = 2 n - (n + 1) = n - 1
  \label{lihsb}
\end{equation}
The usual defintion of $\tmop{Li}_n (t)$ requires analytic continuation at $t
= 1$ but this is not necessary because the hypergeometric function converges
absolutely on the unit circle when it is at least $1$-balanced (\ref{balance})
which is true $\forall n \geqslant 2$. The only \tmtextit{Saalsch\"utzian}
polylogarithm is $\tmop{Li}_2 (t$) {\cite[Eq3.8]{geometricFproperties}}
{\cite[25:12]{functionatlas}}{\cite[1.4.2]{polylogStructure}}

\subsubsection{The Differential Equation Solved by $\tmop{Li}^F_n (t)$ and
Some Combinatorics}

Some combinatorial functions need to be defined before writing the
differential equation solved by $\tmop{Li}_n^F (t)$. Let a {\tmem{partition}}
be an arrangement of the set of elements $1, \ldots, k$ into $n$ subsets where
each element is placed into exactly one set. The number of partitions of the
set $1, \ldots, k$ into $n$ subsets is given by the {\tmem{Stirling numbers of
the second kind}} {\cite[1.1.3]{combPerm}}{\cite[2.7]{ica}} defined by
\begin{equation}
  \begin{array}{ll}
    \left\{ \begin{array}{l}
      k\\
      n
    \end{array} \right\}_{} & = \sum_{j = 0}^n \frac{j^k}{n! (- 1)^{j - n}} 
    \left( \begin{array}{l}
      n\\
      j
    \end{array} \right) \\
    & = \sum_{j = 0}^n \frac{j^k (- 1)^{n - j} }{\Gamma (j + 1) \Gamma (n - j
    + 1)}\\
    & = \frac{(- 1)^{n + 1}}{\Gamma (n)} _k F_{k - 1} \left( \begin{array}{l}
      1 - n, \vec{2}_{k - 1}\\
      \vec{1}_{k - 1}
    \end{array} \right)
  \end{array}
\end{equation}
The $_k F_{k - 1}$ representation of $\left\{ \begin{array}{l}
  k\\
  n
\end{array} \right\}$ is ($n - k$)-balanced (\ref{balance}) since $(k - 1) -
((1 - n) + 2 (k - 1)) = n - k$. The {\tmem{$r$-restricted Stirling numbers of
the second kind}} $\left\{ \begin{array}{l}
  k\\
  n
\end{array} \right\}_r$, or simply the {\tmem{$r$-Stirling numbers}}, counts
the number of partitions of the set 1,...,n into $k$ subsets with the
restriction that the numbers $1, \ldots, r$ belong to distinct subsets.
{\cite{nrs}} The recursion satisfied by $\left\{ \begin{array}{l}
  k\\
  n
\end{array} \right\}_r$ is given by
\begin{equation}
  \begin{array}{ll}
    \left\{ \begin{array}{l}
      k\\
      n
    \end{array} \right\}_r & = \left\{ \begin{array}{ll}
      0 & k < r\\
      \delta_{n, r} & k = r\\
      n \left\{ \begin{array}{l}
        k - 1\\
        n
      \end{array} \right\}_r + \left\{ \begin{array}{l}
        k - 1\\
        n - 1
      \end{array} \right\}_r & n > r
    \end{array} \right.
  \end{array}
\end{equation}
where $\delta_{n, m} = \left\{ \begin{array}{ll}
  1 & n = m\\
  0 & n \neq m
\end{array} \right.$ is the Kronecker delta. Specifically, the
{\tmem{$2$-restricted Stirling numbers}}{\cite[A143494]{oeis}} appearing in
the differential equation for $\tmop{Li}_n^F (t)$ are given by
\begin{equation}
  \begin{array}{lll}
    \left\{ \begin{array}{l}
      k\\
      n
    \end{array} \right\}_2 &  & = \left\{ \begin{array}{l}
      k\\
      n
    \end{array} \right\} - \left\{ \begin{array}{l}
      k - 1\\
      n
    \end{array} \right\}\\
    &  & = \frac{1}{(k - 2) !} \sum_{j = 0}^{k - 2} (- 1)^{j - k}
    \left(\begin{array}{c}
      k - 2\\
      j
    \end{array}\right) (j + 2)^{n - 2}\\
    &  & = (- 1)^k \sum_{j = 0}^{k - 2} \frac{(j + 2)^{n - 2} (- 1)^j}{j! (k
    - 2 - j) !}
  \end{array}
\end{equation}
The ($n + 1$)-th order hypergeometric differential equation (\ref{hode})
satisfied by  f(t)=$\tmop{Li}^F_n (t)$ (\ref{polylog}) \
\begin{equation}
  0 = \left\{ \begin{array}{ll}
    f (t) + \frac{\mathd}{\mathd t} f (t) (t^2 - t) & n = 0\\
    \frac{\mathd}{\mathd t} f (t) + \sum_{m = 2}^{n + 1}  \left(
    \frac{\mathd}{\mathd t^m} f (t) \right)  \left( t^{m - 1} \left\{
    \begin{array}{l}
      n + 1\\
      m
    \end{array} \right\} - t^{m - 2} \left\{ \begin{array}{l}
      n + 1\\
      m
    \end{array} \right\}_2  \right) & n \geqslant 1
  \end{array} \right. \label{polyhde}
\end{equation}
has a most general solution of the form

\begin{equation}
  \begin{array}{l}
    \begin{array}{l}
      f \left( t \right) = x + y_{} G_n (t) + \sum_{m = 1}^{n - 1} z_m \ln
      (t)^m
    \end{array} \label{plhodesol}
  \end{array}
\end{equation}

where $x, y, z_1, \ldots, z_{n - 1}$ are arbitrary parameters and $G_n (t)$
satifies the recursion
\begin{equation}
  \begin{array}{ll}
    G_n (t) & = \left\{ \begin{array}{ll}
      \frac{t}{1 - t} & n = 0\\
      \ln (t - 1) & n = 1\\
      \tmop{Li}_2 (1 - t) + \ln (t - 1) \ln (t) & n = 2\\
      \int \frac{G_{n - 1} (t)}{t} \mathd t & n \geqslant 3
    \end{array} \right.
  \end{array} \label{phi}
\end{equation}
which has the explicit solution
\begin{equation}
  \begin{array}{ll}
    G_n (t) & = \left\{ \begin{array}{ll}
      \frac{t}{1 - t} & n = 0\\
      \ln (t - 1) & n = 1\\
      \tmop{Li}_2 (1 - t) + \ln (t - 1) \ln (t) & n = 2\\
      \frac{(\ln (t - 1) - \ln (1 - t)) \ln (t)^{n - 1}}{\Gamma (n)} +
      \frac{\pi^2}{6} \frac{\ln (t)^{n - 2} }{\Gamma (n - 1)} - \tmop{Li}_n
      (t) & n > 2
    \end{array} \right.
  \end{array}
\end{equation}
The indicial equation of (\ref{polyhde}) at the $t = 1$ is
\begin{equation}
  \begin{array}{l}
    \tmop{ind} (\tmop{Li}_n^F (t)) = - \frac{t (- 1)^{n - 1} \Gamma (n - 1 -
    t) (t - n + 1)^2}{\Gamma (1 - t)}
  \end{array} \label{polyind}
\end{equation}
The ($n + 1$) roots of $\tmop{ind} (\tmop{Li}_n^F (t))$ are the exponents of
(\ref{polyhde}) which are simply
\begin{equation}
  \left\{ t : \tmop{ind} (\tmop{Li}_n^F (t)) = 0\}= 0, 1, \ldots, n - 1, n - 1
  \right.
\end{equation}
where the last root $n - 1$ of $\tmop{ind} (\tmop{Li}_n^F (t))$ is the balance
of $\tmop{Li}_n^F (t)$ (\ref{lihsb}) having multiplicity 2 thus inducing the
logarithmic terms of (\ref{plhodesol}). {\cite[15.31 and 16.33]{ode}} These
equations were derived by writing the differential equation for increasing
values of $n$ and then noticing that the developing pattern of coefficients
were combinatorial. After deriving the general combinatorial differential
equation, it was solved for increasing values of $n$ which resulted in nested
integrals of prior solutions and then the general solution was derived from
that pattern.

\subsubsection{The ``Hypergeometric Form'' of the Zeta Function}

The main focus will be on the special case $\tmop{Li}^F_n (t)$ at unit
argument where it coincides with the Riemann Zeta function at the integers. As
with $\tmop{Li}^F_n (t)$, the symbol $\zeta^F (n)$ refers specifically to the
hypergeometric representation of $\zeta (n)$ at non-negative integer values of
$n$. Using (\ref{hurwitzriemann}) and (\ref{polylog}), it can easily be seen
that $\zeta (n$) can be expressed as a generalized hypergeometric function
(\ref{hypergeom}) with
\begin{equation}
  \begin{array}{ll}
    \zeta^F (n) & = \tmop{Li}^F_n (1)\\
    & ={_{n + 1} F_n} \left( \begin{array}{l}
      \vec{1}_{n + 1}\\
      \vec{2}_n
    \end{array}   \right)\\
    & = \sum_{k = 0}^{\infty} \frac{1}{k!} \frac{\prod_{i = 1}^{n + 1}
    (1)_k}{\prod_{j = 1}^n (2)_k}\\
    & = \sum_{k = 0}^{\infty} \frac{(1)^{n + 1}_k}{k! (2)^n_k}\\
    & = \sum_{k = 0}^{\infty} \frac{\Gamma (k + 1)^n}{\Gamma (k + 2)^n}\\
    & = \sum_{k = 0}^{\infty} (k + 1)^{- n}\\
    & = \zeta (n, 1)\\
    & = \zeta (n)
  \end{array} \label{hyperzeta}
\end{equation}
The value $\zeta^F (0) =_1 F_0 \left( \begin{array}{l}
  1\\
  
\end{array} | 1 \right)$ is singular and so must be calculated with the
reflection equation (\ref{lireflect}) to get $\tmop{Li}^F_0 (- 1) =_1 F_0
\left( \begin{array}{l}
  1\\
  
\end{array} | - 1 \right) = - \frac{1}{2} = \zeta (0)$ which agrees with the
integral form of $\zeta (t) \forall t \neq 1$
\begin{equation}
  \begin{array}{lllll}
    \left. \zeta (t) \right|_{t = 0} & \left. = \left( \frac{1}{2} +
    \frac{1}{t - 1} + 2 \int_0^{\infty} \frac{\sin (s \arctan (s)) (1 +
    s^2)^{- \frac{s}{2}}}{e^{2 \pi s} - 1} \mathd s \right) \right|_{t = 0} &
    = & \tmop{Li}_0 (- 1) = - \frac{1}{2} & 
  \end{array}
\end{equation}

\section{Number Theory, Continued Fractions, and Fractal Strings}

\subsection{Fractal Strings and Dynamical Zeta Functions}

A fractal string is defined as a nonempty open subset of the real line $\Omega
\subseteq \mathbbm{R}$ which can be expressed as a disjoint union of open
intervals $I_j$ being the connected components of $\Omega$.
{\cite[3.1]{fractalzetastrings}}{\cite{weylberry}}{\cite{fdisp}}{\cite{hearfractaldrumshape}}{\cite{gmc}}{\cite{ncfg}}
\begin{equation}
  \begin{array}{ll}
    \Omega & = \bigcup_{j = 1}^{\infty} I_j
  \end{array}
\end{equation}
The length of the $j$-th interval $I_j$ will be denoted $\ell_j$. It will be
assumed that $\Omega$ is standard, meaning that its length is finite, and that
$\ell_j$ is a nonnegative monotically decreasing sequence.
\begin{equation}
  \begin{array}{l}
    | \Omega |_d = \sum_{j = 1}^{\infty} (\ell_j^{})^d < \infty \exists d >
    0\\
    \ell_1 \geqslant \ell_2 \geqslant \ldots \geqslant \text{$\ell_j \geqslant
    \ell_{j + 1} \geqslant \cdots \geqslant 0$}
  \end{array}
\end{equation}
where $\exists d > 0$ means there is at least one value of $d$ for which the
statement is true. It can be the case that $\ell_j = 0$ for some $j$ in which
case $\ell_j$ is a finite sequence. The sequence of lengths of the components
of the fractal string is denoted by
\begin{equation}
  \begin{array}{ll}
    \mathcal{L} & =\{\ell_j \}_{j = 1}^{\infty}
  \end{array}
\end{equation}
The boundary of $\Omega$ in $\mathbbm{R}$ will be denoted by $\partial \Omega
= K \subset \Omega$ which will also denote the boundary of $\mathcal{L}$. Any
totally disconnected bounded perfect subset $K \subset \mathbbm{R}$, or
generally, any compact subset $K \subset \mathbbm{R}$, can be represented as a
string of finite length $| \Omega |_1$. A subset $K$ of a topological space
$\Omega$ is said to be perfect if it is closed and each of its points is a
limit point. Since here $\Omega$ is a metric space and $K \subset \Omega$ is
closed, the Cantor-Bendixon lemma states that there exists a perfect set $P
\subset K$ such that $K - P$ is a most countable. {\cite[2.2 Ex17]{mti}} As
such, $\Omega$ can be defined as the complenent of $K$ in its closed convex
hull, that is, $\Omega = \Omega (K)$ is the smallest compact interval $[a, b]$
containing $K$. The connected components of the bounded open set $\Omega = (a,
b) \backslash K$ are the intervals $I_j$ of the fractal string $\mathcal{L}$
associated with $K$.

\subsection{The Gauss Map $h (x)$}

\begin{figure}[h]
  \resizebox{12cm}{12cm}{\includegraphics{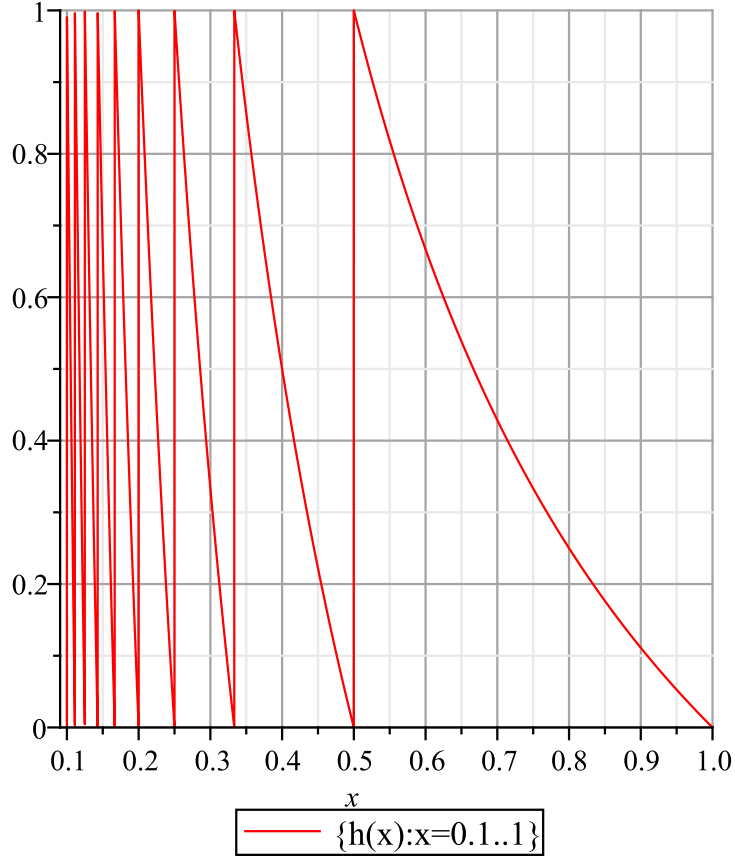}}
  \caption{\label{gaussmapfig}The Gauss Map}
\end{figure}

Let $\Omega_h = (0, 1) \backslash \partial \Omega_h$ where $\partial \Omega_h
= \left\{ \pm \infty, 0, \frac{1}{n} : n \in \mathbbm{Z} \right\}$ is the set
of discontinous boundary points of the Gauss map $h (x) \in \Omega_h \forall x
\nin \partial \Omega_h$, also known as the Gauss function or Gauss
transformation, which maps unit intervals onto unit intervals and by iteration
gives the continued fraction expansion of a real number
\begin{equation}
  \begin{array}{ll}
    h (x) & = \frac{1}{x} - \left\lfloor \frac{1}{x} \right\rfloor\\
    & = - \frac{\left\lfloor \frac{1}{x} \right\rfloor x - 1}{x}\\
    & =\{x^{- 1} \}\\
    & = \frac{1}{x} \tmop{mod} 1
  \end{array}
\end{equation}
where $\left\lfloor x \right\rfloor$ is the floor function, the greatest
integer$\leqslant x$ and $\{x\}= x - \left\lfloor x \right\rfloor$ is the
fractional part of $x$.
{\cite[2.1,3.9.1,9.1,9.3,9.7.1]{cet}}{\cite[A.1.7]{Ctsa}} Clearly $h (x)$ is
also defined outside of $\Omega$
\begin{equation}
  h (x) = \left\{ \begin{array}{ll}
    \frac{1}{x} & x > 1\\
    h (x) & - 1 \leqslant x \leqslant 1\\
    \frac{1}{x} + 1 & x < - 1
  \end{array} \right.
\end{equation}
since
\begin{equation}
  \begin{array}{l}
    
  \end{array} \left\lfloor \frac{1}{x} \right\rfloor = \left\{
  \begin{array}{ll}
    0 & x > 1\\
    \left\lfloor \frac{1}{x} \right\rfloor & - 1 \leqslant x \leqslant 1\\
    - 1 & x < - 1
  \end{array} \right.
\end{equation}
As can be seen in Figure \ref{gaussmapfig}, $h (x)$ is discontinuous at a
countably infinite set of points of Lebesgue measure zero on its boundary
$\partial \Omega_h$
\begin{equation}
  \text{$\left\{ y : \lim_{x \rightarrow y^-} h (x) \neq \lim_{x \rightarrow
  y^+} h (x) \right\} =$} \partial \Omega_h = \left\{ \pm \infty, 0,
  \frac{1}{n} : n \in \mathbbm{Z} \right\}
\end{equation}
The left and right limits of $h (x)$ when $x$ approaches an element on the
boundary $\partial \Omega_h$ is given by
\begin{equation}
  \begin{array}{lll}
    \lim_{x \rightarrow \partial \Omega^-} h (x) & = 0 & \\
    \lim_{x \rightarrow \partial \Omega^{+_{}}} h (x) & = 1 & 
  \end{array}
\end{equation}

\subsubsection{The Frobenius-Perron Transfer Operator}

The Frobenius-Perron transfer operator
{\cite[2.4.4]{Ctsa}}{\cite[2.3.3]{cds}}{\cite[1.3.1,8.2]{cet}}{\cite[1.8,2.4]{dzf}}
of a unit interval mapping $f (y)$ describes how a probablility density $\rho
(y)$ transforms under the action of the map.
\begin{equation}
  \begin{array}{ll}
    {}[U_f \rho] (x) & = \int \delta (x - f (y)) \rho (y) \mathd y
  \end{array}
\end{equation}
where $\delta$ is the Dirac delta function and $\theta$ is the Heaviside step
function.
\begin{equation}
  \begin{array}{ll}
    \text{$\int \delta (x) \mathd x$} & = \theta (x)\\
    \theta (x) & = \left\{ \begin{array}{ll}
      0 & x < 0\\
      1 & x \geqslant 0
    \end{array} \right.
  \end{array} \label{deltastep}
\end{equation}
The function $f (y)$ is the map being iterated and $\rho (y)$ is some density
on which the transfer operator $U$ acts. Essentially, iteration of the map
transforms points to points and iteration of the transfer operator maps point
densities to point densities. The Gauss-Kuzmin-Wirsing(GKW) operator is
obtained by applying the transfer operator to the Guass map.
{\cite[2]{gkwzeta}} {\cite{gklw}} {\cite{gkwoperator}}
\begin{equation}
  \begin{array}{ll}
    {}[U_h \rho] (x) & = \sum_{n = 1}^{\infty} \frac{\rho \left( \frac{1}{n +
    x^{}} \right)}{(n + x)^2}
  \end{array} \label{gkw}
\end{equation}
By changing the variables and order of integration in (\ref{gaussmap}) an
operator equation for $\zeta (s$) is obtained.
\begin{equation}
  \begin{array}{ll}
    \zeta (s) & = \frac{s}{s - 1} - s \int_0^1 x [U_h x^{s - 1}] \mathd x\\
    & = \frac{s}{s - 1} - s \int_0^1 x \int \delta (x - h (y)) y^{s - 1}
    \mathd y \mathd x\\
    & = \frac{s}{s - 1} - s \int_0^1 x \int \delta (x - \left( y^{- 1} -
    \left\lfloor y^{- 1} \right\rfloor \right)) y^{s - 1} \mathd y \mathd x\\
    & = \frac{s}{s - 1} - s \int_0^1 x \sum_{n = 1}^{\infty} \frac{\left(
    \frac{1}{n + x^{}} \right)^{s - 1}}{(n + x)^2} \mathd x
  \end{array}
\end{equation}
An operator similiar to (\ref{gkw}) is
\begin{equation}
  \begin{array}{ll}
    {}[S_h \rho] (x) & = \sum_{n = 1}^{\infty} \rho \left( \frac{1}{n} \right)
    - \rho \left( \frac{1}{n + x} \right)
  \end{array} \label{S}
\end{equation}
The action of $[U_h \rho] (x)$ on the identity function $x \rightarrow x$ is
given by
\begin{equation}
  \begin{array}{ll}
    {}[U_h x] (x) & = \sum_{n = 1}^{\infty} \frac{\frac{1}{n + x^{}}}{(n +
    x)^2}\\
    & = \sum_{n = 1}^{\infty} \frac{1}{(n + x)^3}\\
    & = - \frac{\Psi^{(2)} (x + 1)}{2}\\
    & = \frac{_4 F_3 \left( \begin{array}{llll}
      1 & x + 1 & x + 1 & x + 1\\
      & x + 2 & x + 2 & x + 2
    \end{array} \right)}{(x + 1)^3} 
  \end{array}
\end{equation}
where $\Psi^{(n)} (x)$ is the polygamma function (\ref{phiinf}). The area
under the curve of $[U_h x] (x)$ is
\begin{equation}
  \begin{array}{l}
    \int_0^1
  \end{array} [U_h x] (x) \mathd x = \begin{array}{l}
    \int_0^1
  \end{array} - \frac{\Psi^{(2)} (x + 1)}{2} \mathd x = \frac{1}{2}
\end{equation}
The identity action of $[S_h \rho] (x)$ is
\begin{equation}
  \begin{array}{ll}
    {}[S_h x] (x) & = \sum_{n = 1}^{\infty} \frac{1}{n} - \frac{1}{n + x}\\
    & = \gamma + \Psi (x + 1)
  \end{array} \label{Seye}
\end{equation}
where $\matheuler$ is Euler's constant
\begin{equation}
  \begin{array}{ll}
    \matheuler & = \lim_{n \rightarrow \infty} \sum_{k = 1}^n \frac{1}{k} -
    \ln (n)\\
    & = \lim_{s \rightarrow 1} \zeta (s) - \frac{1}{s - 1}\\
    & = \lim_{s \rightarrow 1} \zeta (s) + \int_1^{\infty} h (x) x^{s - 1}
    \mathd x\\
    & = \lim_{s \rightarrow 1} \frac{1}{s - 1} - s \int_0^1 h (x) x^{s - 1}
    \mathd x + \int_1^{\infty} h (x) x^{s - 1} \mathd x\\
    & = 1 - \int_0^1 h (x) \mathd x\\
    & \approx 0.57721566490153286060651209
  \end{array}
\end{equation}
and the area under its curve is given by
\begin{equation}
  \begin{array}{ll}
    \int_0^1 [S_h x] (x) \mathd x & = \int_0^1 \gamma + \Psi (x + 1) \mathd x
    = 1 - \gamma
  \end{array} \label{Sarea}
\end{equation}

\subsubsection{Piecewsise Integration of $h (x)$}

The Guass map $h (x) \in \Omega_h$ is piecewise monotone {\cite[2.1]{dzf}}
between the points of $\partial \Omega_h$, and thus partitions the unit
interval infinite covering set of decreasing open intervals seperated by
$\partial \Omega_h$. {\cite[5.7.1]{cet}} Let $I_n$ be an infinite set of open
intervals
\begin{equation}
  \begin{array}{ll}
    I_n & = \left\{ \begin{array}{ll}
      (1, \infty) & n = 0\\
      \left( \frac{1}{n + 1}, \frac{1}{n} \right) & 0 < n < \infty\\
      (0, 0) = \emptyset & n = \infty
    \end{array} \right.
  \end{array}
\end{equation}
It is easy to see that
\begin{equation}
  \begin{array}{lll}
    \Omega_h \cup \partial \Omega_h = & [0, 1] & = \text{$\bigcup_{n =
    1}^{\infty} I_n$}\\
    & [0, \infty] & = \text{$\bigcup_{n = 0}^{\infty} I_n$}
  \end{array} \label{unitpart}
\end{equation}
Define the Gauss map partition $h_n (x)$ where $\{h_n (x) \neq 0 : x \in I_n
\}$ as a piecewise step function
\begin{equation}
  \begin{array}{ll}
    h_n (x) & = \left\{ \begin{array}{ll}
      \frac{1 - x n}{x} & \frac{1}{n + 1} < x < \frac{1}{n}\\
      0 & \tmop{otherwise}
    \end{array} \right.\\
    & = \left( \frac{1 - x n}{x} \right) \left( \theta \left( \frac{x n + x -
    1}{n + 1} \right) - \theta \left( \frac{x n - 1}{n} \right) \right)
  \end{array} \label{hn}
\end{equation}
where $\theta (t)$ is the Heaviside step function (\ref{deltastep}). We can
reassemble all of the $\{h_n (x)\}_{n = 1}^{\infty}$ to recover $h (x$)
\begin{equation}
  \begin{array}{ll}
    h (x) & = \sum_{n = 1}^{\infty} h_n (x)\\
    & = \sum_{n = 1}^{\infty} \left( \frac{1 - x n}{x} \right) \left( \theta
    \left( \frac{x n + x - 1}{n + 1} \right) - \theta \left( \frac{x n - 1}{n}
    \right) \right)
  \end{array} \label{hpw}
\end{equation}
where only one of the $h_n (x)$ is $\tmop{nonzero}$ for each $x$. By setting
$n = \left\lfloor \frac{1}{x} \right\rfloor$ in (\ref{hn}) we get
\begin{equation}
  \begin{array}{ll}
    h (x) & = \left\{ \begin{array}{ll}
      \frac{1 - x \left\lfloor \frac{1}{x} \right\rfloor}{x} &
      \frac{1}{\left\lfloor \frac{1}{x} \right\rfloor + 1} < x <
      \frac{1}{\left\lfloor \frac{1}{x} \right\rfloor}\\
      0 & \tmop{otherwise}
    \end{array} \right.\\
    & = \left( \frac{1 - x \left\lfloor \frac{1}{x} \right\rfloor}{x} \right)
    \left( \theta \left( \frac{x \left\lfloor \frac{1}{x} \right\rfloor + x -
    1}{\left\lfloor \frac{1}{x} \right\rfloor + 1} \right) - \theta \left(
    \frac{x \left\lfloor \frac{1}{x} \right\rfloor - 1}{\left\lfloor
    \frac{1}{x} \right\rfloor} \right) \right)
  \end{array}
\end{equation}
Define the partitioned integral operator $\left[ P f (x) ; x \right] (n)$ by
\begin{equation}
  \begin{array}{lllll}
    \text{$\left[ P_{} f (x) ; x \right] (n)$} & = \text{$\left[ P_{} f (x)
    \right] (n)$} & = \left[ P f \right] (n) & = \int_{\frac{1}{n +
    1}}^{\frac{1}{n}} f (x) \mathd x & = \int_{I_n}^{} f (x) \mathd x
  \end{array}
\end{equation}
where by convention we have
\begin{equation}
  \begin{array}{llll}
    \left[ P f (x) \right] (0) & = \lim_{n \rightarrow 0^+} \int_{I_n}^{} f
    (x) \mathd x & = \int_1^{\infty} f (x) \mathd x & \\
    \left[ P_{} f (x) \right] (\infty) & = \lim_{n \rightarrow \infty}
    \int_{I_n}^{} f (x) \mathd x & = \int_0^0 f (x) \mathd x & = 0
  \end{array}
\end{equation}
Thus
\begin{equation}
  \begin{array}{lll}
    \int_0^1 f (x) \mathd x & = \sum_{n = 1}^{\infty} \left[ P_{} f (x)
    \right] (n) & = \sum_{n = 1}^{\infty} \int_{I_n}^{} f (x) \mathd x\\
    \int_0^{\infty} f (x) \mathd x & = \sum_{n = 0}^{\infty} \left[ P_{} f (x)
    \right] (n) & = \sum_{n = 0}^{\infty} \int_{I_n}^{} f (x) \mathd x
  \end{array}
\end{equation}
Each interval $I_n$ has the length
\begin{equation}
  \begin{array}{ll}
    \ell_{} I_n & = \left[ P 1 \right] (n)\\
    & = \int_{\frac{1}{n + 1}}^{\frac{1}{n}} 1 \mathd x\\
    & = \frac{1}{n} - \frac{1}{n + 1}\\
    & = \frac{1}{n (n + 1)}
  \end{array} \label{obblongarea}
\end{equation}
The elements $n (n + 1)$ are known as the oblong numbers
{\cite[A002378]{oeis}}. It is seen, together with (\ref{Seye}), that
\begin{equation}
  \begin{array}{ll}
    | \Omega_h |_{} & = \int_0^1 1 \mathd x\\
    & = \sum_{n = 1}^{\infty} \ell_{} I_n\\
    & = \sum_{n = 1}^{\infty} \frac{1}{n (n + 1)}\\
    & = [S_h x] (1)\\
    & = \gamma + \Psi (2)\\
    & = 1
  \end{array} \label{hcover}
\end{equation}
The piecewise integral operator $\left[ P_{} f (x) ; x \right] (n)$ can be
used to calculate the area under the curve of $h (x)$ which is also equal to
the area under the curve of $[S_h x] (x)$. Let the length of the $n$-th
component $h_n (x)$ be denoted by
\begin{equation}
  \begin{array}{ll}
    \ell h_n & = \left[ P_{} h (x) ; x \right] (n)\\
    & = \int^{}_{I_n} h (x) \mathd x\\
    & = \int_{I_n}^{} h_n (x) \mathd x\\
    & = \int_0^1 h_n (x) \mathd x\\
    & = \frac{\ln (n + 1) n + \ln (n + 1) - \ln (n) n - \ln (n) - 1}{n + 1}
  \end{array}
\end{equation}
Regarding $h (x)$ as a fractal string $\mathcal{L}_h =\{h_n (x)\}_{n =
1}^{\infty}$ its length $|\mathcal{L}_h |$ is given by
\begin{equation}
  \begin{array}{ll}
    \text{$|\mathcal{L}_h |$} & = \int_0^1 h (x) \mathd x\\
    & = \sum_{n = 1}^{\infty} \ell h_n\\
    & = \sum_{n = 1}^{\infty} \frac{\ln (n + 1) n + \ln (n + 1) - \ln (n) n -
    \ln (n) - 1}{n + 1}\\
    & = \int_0^1 [S_h x] (x) \mathd x\\
    & = \int_0^1 \gamma + \Psi (x + 1) \mathd x\\
    & = 1 - \gamma
  \end{array}
\end{equation}
If $n = 0$ in (\ref{unitpart}) we get the interval $\left. I_0 = \left(
\frac{1}{1}, \frac{1}{0} \right) = (1, \infty \right)$ and
\begin{equation}
  \begin{array}{ll}
    \ell h_0 = & \int_{I_0}^{} h (x) \mathd x = \int_1^{\infty} \frac{1}{x}
    \mathd x = \infty
  \end{array}
\end{equation}
but if we choose a finite cutoff then
\begin{equation}
  \begin{array}{ll}
    \int_1^y h (x) \mathd x & = \int_1^y \frac{1}{x} \mathd x\\
    & = \ln (y)
  \end{array}
\end{equation}
and
\begin{equation}
  \begin{array}{l}
    \int_1^{\infty}
  \end{array} \frac{\ln (y)}{y^n} \mathd y = \frac{1}{(n - 1)^2}
\end{equation}
thus
\begin{equation}
  \begin{array}{ll}
    \sum_{n = 2}^{\infty} \begin{array}{l}
      \int_1^{\infty}
    \end{array} \frac{\int_1^y h (x) \mathd x}{y^n} \mathd y & = \sum_{n =
    2}^{\infty} \begin{array}{l}
      \int_1^{\infty}
    \end{array} \frac{\int_1^y \frac{1}{x} \mathd x}{y^n} \mathd y\\
    & = \sum_{n = 2}^{\infty} \begin{array}{l}
      \int_1^{\infty}
    \end{array} \frac{\ln (y)}{y^n} \mathd y\\
    & = \sum_{n = 2}^{\infty} \frac{1}{(n - 1)^2}\\
    & = \zeta (2)\\
    & = \frac{\pi^2}{6}
  \end{array}
\end{equation}

\subsubsection{The Mellin Transform}

The Mellin transform {\cite[II.10.8]{mmp}}{\cite[3.6]{piagm}} is defined as
\begin{equation}
  \begin{array}{ll}
    M^{(a, b)}_{x \rightarrow s} f (x) & = \int_a^b f (x) x^{s - 1} \mathd x
  \end{array} \label{mellin}
\end{equation}
where the usual definition of the Mellin transform is $M^{(0, \infty)}_{x
\rightarrow s} f (x)$. Somewhat incredibly, by taking the Mellin
transformation of $h (x)$ over the unit interval, we get an analytic
continuation of $\zeta (s)$ which is convergent when $s$ is not equal to a
negative integer, $0$, or $1$. When $s$ is a negative integer or 0 the limit
or analytic continuation must be taken since the series is formally divergent
at these points, and of course the series $s = 1$ diverges.
{\cite{newtonzeta}} {\cite{gkwzeta}} {\cite{yarh}}
\begin{equation}
  \begin{array}{ll}
    M_{x \rightarrow s}^{I_n} h (x) & = \left[ P_{} h (x) x^{s - 1} ; x
    \right] (n)\\
    & = \int_{\frac{1}{n + 1}}^{\frac{1}{n}} \left( \frac{1}{x} -
    \left\lfloor \frac{1}{x} \right\rfloor \right) x^{s - 1} \mathd x\\
    & = - \frac{n (n + 1)^{- s} + s (n + 1)^{- s} - n^{1 - s}}{s (s - 1) }
  \end{array}
\end{equation}
\begin{equation}
  \begin{array}{ll}
    \zeta (s) & = \frac{1}{s - 1} - s M_{x \rightarrow s}^{(0, 1)} h (x)\\
    & = \frac{1}{s - 1} - s \int_0^1 h (x) x^{s - 1} \mathd x\\
    & = \frac{1}{s - 1} - s \int_0^1 \left( \frac{1}{x} - \left\lfloor
    \frac{1}{x} \right\rfloor \right) x^{s - 1} \mathd x\\
    & = \frac{s}{s - 1} - s \sum_{n = 1}^{\infty} M_{x \rightarrow s}^{I_n} h
    (x)\\
    & = \frac{s}{s - 1} - s \sum_{n = 1}^{\infty} - \frac{n (n + 1)^{- s} + s
    (n + 1)^{- s} - n^{1 - s}}{s (s - 1) }
  \end{array} \label{gaussmap}
\end{equation}
The term $\frac{1}{s - 1}$ changes to $\frac{s}{s - 1} = \text{$\frac{1}{s -
1} - (- 1)$}$ by subtracting the residue
{\cite[10.41]{rca}}{\cite[6.1]{whittaker-watson}} of
\begin{equation}
  M_{x \rightarrow s}^{I_0} h (x) = \int_{I_0}^{} h (x) x^{s - 1} \mathd x =
  \int_1^{\infty} \frac{x^{s - 1}}{x} \mathd x = - \frac{1}{s - 1}
\end{equation}
at the singular point $s = 1$, which happens to coincide with $\sum_{s =
2}^{\infty} \frac{- \frac{1}{s - 1}}{s}$
\begin{equation}
  \begin{array}{ll}
    \tmop{Res} \left( \int_1^{\infty} \frac{x^{s - 1}}{x} \mathd x ; 1 \right)
    & = \tmop{Res} \left( - \frac{1}{s - 1} ; 1 \right)\\
    & = \sum_{s = 2}^{\infty} - \frac{\frac{1}{s - 1} }{s}\\
    & = - 1
  \end{array}
\end{equation}

\subsection{The Harmonic Sawtooth w(x)}

Define the harmonic sawtooth map $w (x) \in \Omega_h \backslash \partial
\Omega_h$ which shares the same domain and boundary as the Gauss map $h (x)$
to which it is similiar, and also has the property that its Mellin transform
is the (appropriately scaled) zeta function. The $n$-th component $w_n (x)$ is
defined over the $n$-th interval $I_n$
\begin{equation}
  \begin{array}{ll}
    w_n (x) & = \left\{ \begin{array}{ll}
      n (x n + x - 1) & \frac{1}{n + 1} < x < \frac{1}{n}\\
      0 & \tmop{otherwise}
    \end{array} \right.\\
    & = n (x n + x - 1) \left( \theta \left( \frac{x n + x - 1}{n + 1}
    \right) - \theta \left( \frac{x n - 1}{n} \right) \right)
  \end{array}
\end{equation}
and by the substitution $n \rightarrow \left\lfloor \frac{1}{x} \right\rfloor$
we have
\begin{equation}
  \begin{array}{ll}
    w (x) & = \sum_{n = 1}^{\infty} w_n (x)\\
    & = \sum_{n = 1}^{\infty} n (x n + x - 1) \left( \theta \left( \frac{x n
    + x - 1}{n + 1} \right) - \theta \left( \frac{x n - 1}{n} \right)
    \right)\\
    & = \left\lfloor \frac{1}{x} \right\rfloor \left( x \left\lfloor
    \frac{1}{x} \right\rfloor + x - 1 \right)
  \end{array}
\end{equation}
Unlike $h (x)$ which is nonzero outside of $|x| > 1$, the (harmonic) sawtooth
map has $w (x) = 0 \forall |x| > 1$.

\begin{figure}[h]
  \resizebox{12cm}{12cm}{\includegraphics{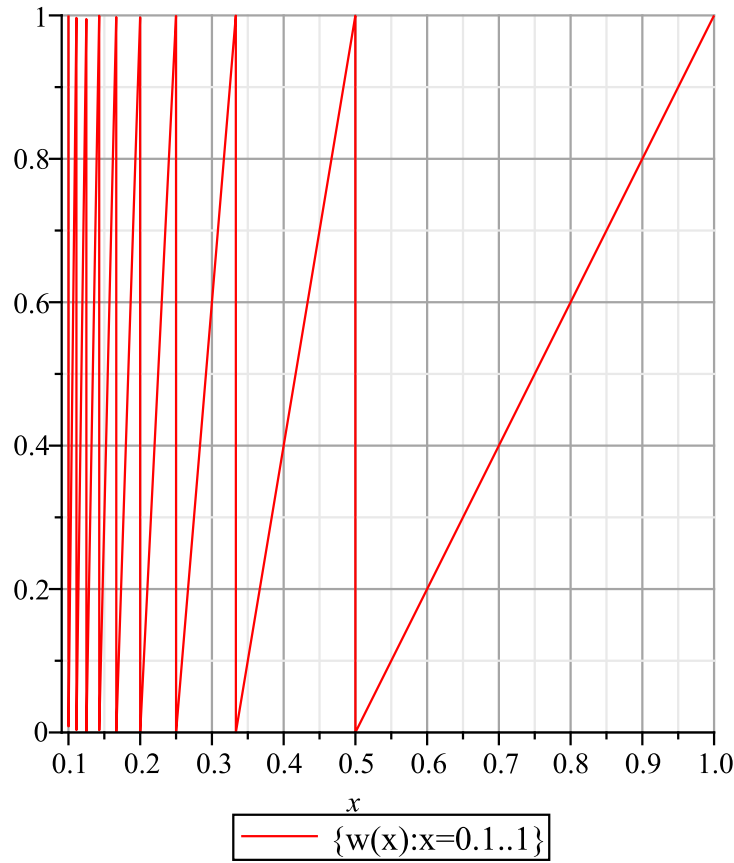}}
  \caption{The Harmonic Sawtooth}
\end{figure}

The length of each component of $w (x)$ is
\begin{equation}
  \begin{array}{ll}
    \ell w_n & = \left[ P_{} w (x) ; x \right] (n)\\
    & = \int^{}_{I_n} w (x) \mathd x\\
    & = \frac{1}{2 (n + 1) n}
  \end{array}
\end{equation}
So that the total length of the harmonic sawtooth string $\mathcal{L}_w$ is
\begin{equation}
  \begin{array}{ll}
    |\mathcal{L}_w | & = \int_0^1 w (x) \mathd x\\
    & = \sum_{n = 1}^{\infty} \ell w_n\\
    & = \sum_{n = 1}^{\infty} \frac{1}{2 (n + 1) n}\\
    & = \frac{1}{2}
  \end{array}
\end{equation}
The infinite set of Mellin transforms of $w_n (x)$
\begin{equation}
  \begin{array}{ll}
    M_{x \rightarrow s}^{I_n} w (x) & = M_{x \rightarrow s}^{(0, 1)} w_n (x)\\
    & = \left[ P_{} w (x) x^{s - 1} ; x \right] (n)\\
    & = \int_{\frac{1}{n + 1}}^{\frac{1}{n}} n (x n + x - 1) x^{s - 1} \mathd
    x\\
    & = \int_0^1 n (x n + x - 1) \left( \theta \left( \frac{x n + x - 1}{n +
    1} \right) - \theta \left( \frac{x n - 1}{n} \right) \right) x^{s - 1}
    \mathd x\\
    & = - \frac{n (n + 1)^{- s} + s (n + 1)^{- s} - n^{1 - s}}{s (s - 1) }
  \end{array}
\end{equation}
are summed to get an expression for $\{\zeta (s) : \mathfrak{R}(s) \nin
\mathbbm{N}_{0^-} \}$
\begin{equation}
  \begin{array}{ll}
    \zeta (s) & = s \frac{s + 1}{s - 1} \int_0^1 \left\lfloor \frac{1}{x}
    \right\rfloor \left( x \left\lfloor \frac{1}{x} \right\rfloor + x - 1
    \right) x^{s - 1} \mathd x\\
    & = \sum_{n = 1}^{\infty} s \frac{s + 1}{s - 1} M_{x \rightarrow s}^{I_n}
    w (x)\\
    & = \sum_{n = 1}^{\infty} s \frac{s + 1}{s - 1} \int_{\frac{1}{n +
    1}}^{\frac{1}{n}} n (x n + x - 1) x^{s - 1} \mathd x\\
    & = \sum_{n = 1}^{\infty} s \frac{s + 1}{s - 1} \left( - \frac{n (n +
    1)^{- s} + s (n + 1)^{- s} - n^{1 - s}}{s (s - 1) } \right)\\
    & = \sum_{n = 1}^{\infty} \frac{n (n + 1)^{- s} - n^{1 - s} + s n^{-
    s}}{s - 1}
  \end{array} \label{sawzeta}
\end{equation}

\subsection{The Prime Numbers}

Let $\mathbbm{P}=\{2, 3, 5, 7, 11, 13, 17, 19, 23, 29, \ldots\}$ denote the
set of prime numbers and $\mathbbm{N}_1 =\{1, 2, \ldots\}$ and $\mathbbm{N}_0
=\{0, 1, 2, \ldots\}, \mathbbm{N}=\{\ldots, - 2, - 1, 0, 1, 2, \ldots\}$ be
the set of positive, non-negative, and signed integers.

\subsubsection{The Prime Counting Function: $\pi (x)$}

The prime counting function $\pi (x)$ counts the number of primes less than a
given number. It can written as
\begin{equation}
  \begin{array}{ll}
    \pi (x) & = \sum_{p < x}^{p \in \mathbbm{P}} 1
  \end{array} \label{primecounting}
\end{equation}
which is essentially a step function which increases by 1 for each prime.
{\cite[15.11]{gamma}}

\subsubsection{von Mangoldt and Chebyshev's Functions: $\Lambda (x), \theta
(x), \text{$\psi (x)$}$}

Chebyshev's function of the first kind $\theta (x)$ is the sum of the
logarithm of all primes $\leqslant x$
\begin{equation}
  \begin{array}{ll}
    \theta (x) & = \sum_{k = 1}^{\pi (x)} \ln (p_k)\\
    & = \ln \left( \sum_{k = 1}^{\pi (x)} p_k \right)
  \end{array}
\end{equation}
where $p_k \in \mathbbm{P}$ is the $k$-th prime. {\cite[4.4]{edwardszeta}} The
generalization of $\pi (x)$ is the Chebyshev function of the second kind
\begin{equation}
  \begin{array}{ll}
    \psi (x) & = \sum_{p^r \leqslant x}^{\{p \in \mathbbm{P}, r \in
    \mathbbm{N}_1 \}} \ln (p)\\
    & = \sum_{k = 1}^{\left\lfloor \log_2 (x) \right\rfloor} \theta
    (x^{\frac{1}{k}})\\
    & = \ln (\tmop{lcm} (1, 2, 3, \ldots, \left\lfloor x \right\rfloor))\\
    & = \sum^{n \leqslant x}_n \Lambda (n)\\
    & = x - \frac{\ln (1 - x^{- 2})}{2} - \ln (2 \pi) - \sum^{\zeta (\rho) =
    0}_{\rho_{}} \frac{x^{\rho_{}}}{\rho_{}} \forall \mathcal{I}(\rho) \neq 0
  \end{array} \label{cheb2}
\end{equation}
where the first sum ranges over the primes $p \in \mathbbm{P}$ and positive
integers $r$ and the sum over $\rho$ is von Mangoldt's formula where $\rho$
ranges over the non-trivial roots of $\zeta (s)$ in increasing order. The
function $\tmop{lcm} (\ldots .)$ is the least common multiple, and $\Lambda
(x)$ is the von Mangoldt function.
\begin{equation}
  \begin{array}{ll}
    \Lambda (x) & = \left\{ \begin{array}{ll}
      \ln (p) & \{n = p^k : p \in \mathbbm{P}, k \in \mathbbm{N}_1 \}\\
      0 & \tmop{otherwise}
    \end{array} \right.\\
    & = \ln \left( \frac{\tmop{lcm} (1, 2, \ldots, n)}{\tmop{lcm} (1, 2,
    \ldots, n - 1)} \right)
  \end{array}
\end{equation}
$\Lambda (s)$ is related to $\zeta (s)$ by
\begin{equation}
  - \frac{\frac{\mathd}{\mathd s} \zeta (s)}{\zeta (s)} = \sum_{n =
  1}^{\infty} \frac{\Lambda (n)}{n^s} \forall \mathcal{R}(s) > 1
\end{equation}
Chebyshev proved that $\pi (x), \theta (x), \tmop{and} \text{$\psi (x)$}$ \
have the same scaled asymptotic limit.
\begin{equation}
  \lim_{x \rightarrow \infty} \frac{\pi (x)}{\left( \frac{x}{\ln (x)} \right)}
  = \lim_{x \rightarrow \infty} \frac{\psi (x)}{x} = \lim_{x \rightarrow
  \infty} \frac{\theta (x)}{x} = 1
\end{equation}
{\cite[1.3]{ent}} {\cite[I.4]{primedist}}
{\cite[4.3-4.4\&3.1-3.2]{edwardszeta}} {\cite{hlrzdp}} {\cite[15.11]{gamma}}.
Note that {\cite{gamma}} incorrectly defines $\psi (x)$ as $\ln (\gcd
(\ldots)) .$

\section{Analytic Continuation}

\subsection{Continuation of $_{n + 1} F_n$ Near Unit Argument}

The continuation formula for Gauss's hypergeometric function $_2 F_1$ near
unit argument is well known

\begin{equation}
  \begin{array}{ll}
    \frac{\Gamma (a_1) \Gamma (a_2)}{\Gamma (b_1)} _2 F_1 \left(
    \begin{array}{ll}
      a_1 & a_2\\
      b _1 & 
    \end{array} | z \right) & = \sum_{n = 0}^{\infty} \frac{(- 1)^n (1 -
    z)^n}{n!}  \frac{\Gamma (a_1 + n) \Gamma (a_2 + n) \Gamma (s_1 -
    n)}{\Gamma (a_1 + s_1) \Gamma (a_2 + s_1)}\\
    & + (1 - z)^{s_1} \sum_{n = 0}^{\infty} \frac{(- 1)^n (1 - z)^n}{n!} 
    \frac{\Gamma (a_1 + s_1 + n) \Gamma (a_2 + s_1 + n) \Gamma (- s_1 -
    n)}{\Gamma (a_1 + s_1) \Gamma (a_2 + s_1)}
  \end{array}
\end{equation}

where $s_1 = b_1 - a_1 - a_2$ is the balance (\ref{balance}) of $_2 F_1$ which
must not be equal to an integer, that is, $_2 F_1$ cannot be $s_1$-balanced. A
function is said to be $k$-balanced only when $k$ is an integer. When
$\mathcal{R}(s_1) > 0$ the value at $z = 1$ is finite and given by the
Gaussian summation formula
\begin{equation}
  \begin{array}{ll}
    \frac{_2 F_1 \left( \begin{array}{ll}
      a_1 & a_2\\
      b_1 & 
    \end{array} \right)}{\Gamma (b_1)} & = \frac{\Gamma (b_1 - a_1 -
    a_2)}{\Gamma (b_1 - a_1) \Gamma (b_1 - a_2)} \label{gsum}\\
    & = \frac{\Gamma (s_1)}{\Gamma (a_1 + s_1) \Gamma (a_2 + s_1)}
  \end{array}
\end{equation}
It is obvious that $\lim_{t \rightarrow 1} \tmop{Li}_1^F (t) = \lim_{t
\rightarrow 1} \text{$_2 F_1 \left( \begin{array}{ll}
  1 & 1\\
  2 & 
\end{array} | t \right) = \zeta^F (1)$} = \infty$ is $0$-balanced and of
course equal to the divergent harmonic series so the continuation formula does
not apply. However, B\"uhring and Srivastava
{\cite{hfunitwolfgang}}{\cite{hypercont}} generalized this relation to all
$_{n + 1} F_n$ by expanding (\ref{gsum}) as a series then interchanging the
order of summations to derive a recurrence with respect to $n$

\begin{equation}
  \begin{array}{ll}
    _{n + 1} F_n \left( \begin{array}{l}
      a_1, \ldots, a_{n + 1}\\
      b_1, \ldots, b_n
    \end{array} | t \right) & = \frac{\Gamma (b_n) \Gamma (b_{n - 1})}{\Gamma
    (a_{n + 1}) \Gamma (b_n + b_{n - 1} - a_{n + 1})}\\
    & \cdot \sum_{m = 0}^{\infty} \frac{(b_n - a_{n + 1})_m (b_{n - 1} - a_{n
    - 1})_m}{(b_n + b_{n - 1} - a_{n + 1})_m m!} _n F_{n - 1} \left(
    \begin{array}{l}
      a_1, \ldots, a_n\\
      b_1, \ldots, b_{n - 2}, b_{n - 1} + b_n - a_{n + 1} + m
    \end{array} | t \right)
  \end{array} \label{Frecur}
\end{equation}

which is valid $\forall \{\mathcal{R}(a_i) > 0 : 1 \leqslant i \leqslant n +
1\}$. The $m$-th term of the summand in (\ref{Frecur}) is contiguous
(\ref{contig}) to the ($m - 1$)-th and ($m + 1$)-th terms and thus a linear
relationship can always be found between neighboring terms.

\subsection{The Continuation of $\tmop{Li}_n^F (t) \tmop{and} \zeta^F (n)$ via
Contiguous Functions}

There are 4 functions contiguous (\ref{contig}) to $\tmop{Li}_n^F (t)$, only 3
of them are unique, and just 1 of them is interesting. The functions are
obtained by shifting one of the numerator parameters $a_i \pm 1$ or shifting
one of the denominator parameters $b_i \pm 1$. Shifting any of the $a$
parameters or any of the $b$ parameters will suffice since they are all equal
and $_p F_q$ is invariant with respect to the ordering of parameters. Let
$\vec{c}^+_n$ and $\vec{c}^-_n$ denote the parameter vector $\vec{c}_n$ where
one element is shifted up or down by $1$.
\begin{equation}
  \begin{array}{ll}
    \vec{c}^+_n & = \vec{c}_{n - 1}, c + 1 = \underbrace{c, \ldots, c}_{n -
    1}, c + 1\\
    \vec{c}^-_n & = \vec{c}_{n - 1}, c - 1 = \underbrace{c, \ldots, c}_{n -
    1}, c - 1
  \end{array} \label{shift}
\end{equation}
For example, $\vec{4}_3^+ = 4, 4, 5$. Two of the four functions contiguous to
$\tmop{Li}_n^F (t)$ are identical
\begin{equation}
  \begin{array}{ll}
    _{n + 1} F_n \left( \begin{array}{l}
      \vec{1}^+_{n + 1}\\
      \vec{2}_n
    \end{array}   | t \right) t & =_{n + 1} F_n \left( \begin{array}{l}
      \vec{1}^{}_{n + 1}\\
      \vec{2}^-_{^{} n}
    \end{array}   | t \right) t = \tmop{Li}_{n - 1}^F (t)
  \end{array}
\end{equation}
Shifting any $a_i$ up is equivalent to shifting any $b_i$ down, both
operations take $\tmop{Li}_n^F (t)$ to $\tmop{Li}_{n - 1}^F (t)$. Shifting any
$a_i$ down results in the identity function since it puts a $0$ in the
numerator.
\begin{equation}
  _{n + 1} F_n \left( \begin{array}{l}
    \vec{1}^-_{n + 1}\\
    \vec{2}_n
  \end{array}   | t \right) t = t
\end{equation}
Thus, the only interesting function continguous to $\tmop{Li}_n^F (t)$ is
obtained by shifting one of the denominator parameters up. Let this function
be denoted by $\tmop{Li}_n^{F + 1} (t)$
\begin{equation}
  \begin{array}{ll}
    \text{$\tmop{Li}_n^{F + 1} (t)$} =_{n + 1} F_n \left( \begin{array}{l}
      \vec{1}^{}_{n + 1}\\
      \vec{2}^+_{^{} n}
    \end{array}   | t \right) & = \left\{ \begin{array}{ll}
      I_0 \left( 2 \sqrt{t} \right) - \frac{1}{\sqrt{t}} I_1 \left( 2 \sqrt{t}
      \right) & n = 0\\
      \frac{e^t}{t} - \frac{1}{t} - 1 & n = 1\\
      (- 1)^n \left( 1 - \frac{\tmop{Li}_1 (t)}{t} + \sum_{k = 1}^{n - 1} (-
      1)^{k + 1} \tmop{Li}_k (t) \right) & n \geqslant 2
    \end{array} \right.
  \end{array}
\end{equation}
where $I_n (x)$ is a {\tmem{modified Bessel function of the first kind}}
{\cite[65]{sf}} {\cite[6.9.1]{htf1}}
\begin{equation}
  \begin{array}{ll}
    I_n (x) & = \frac{x^n }{\Gamma (n + 1) 2^n} _0 F_1 \left( \begin{array}{l}
      \\
      n + 1
    \end{array} | \frac{x^2}{4} \right)
  \end{array}
\end{equation}
Before applying (\ref{Frecur}), the notation will be simplified by extending
(\ref{shift}) so that repeated shifts can be written more easily
\begin{equation}
  \begin{array}{ll}
    \vec{c}^{+ j}_n & = \vec{c}_{n - 1}, c + j = \underbrace{c, \ldots, c}_{n
    - 1}, c + j\\
    \vec{c}^{- j}_n & = \vec{c}_{n - 1}, c - j = \underbrace{c, \ldots, c}_{n
    - 1}, c - j
  \end{array} \label{multishift}
\end{equation}
where clearly $\vec{c}^+_n = \vec{c}^{+ 1}_n$ and $\vec{c}^-_n = \vec{c}^{-
1}_n$. The goal is to extend $\tmop{Li}_n^{F + 1} (t)$ to all $\tmop{Li}_n^{F
+ m} (t)$ by repeated application of $\vec{c}^{+ 1}_n$. Applying
(\ref{Frecur}) to (\ref{polylog}) gives the continuation of $\tmop{Li}_n^F (t)
\rightarrow \tmop{Li}_{n + 1}^F (t) \forall n \geqslant 1$ by setting $a_{1
\ldots n + 1} = \vec{1}_{n + 1}$ and $b_{1 \ldots n} = \vec{2}_n$ which
results in
\begin{equation}
  \begin{array}{ll}
    \tmop{Li}_n^F (t) & ={_{n + 1} F_n} \left( \begin{array}{l}
      \vec{1}_{n + 1}\\
      \vec{2}_n
    \end{array}   | t \right) t \forall n \geqslant 0\\
    & = t \sum_{m = 0}^{\infty} \left( \frac{_n F_{n - 1} \left(
    \begin{array}{l}
      \vec{1}_n\\
      \vec{2}_{n - 2}, 3 + m
    \end{array} | t \right)}{(m + 1) (m + 2)} \right) \forall n \geqslant 2\\
    & = t \sum_{m = 0}^{\infty} \left( \frac{_n F_{n - 1} \left(
    \begin{array}{l}
      \vec{1}_n\\
      \vec{2}^{+ m + 1}_{n - 1}
    \end{array} | t \right)}{(m + 1) (m + 2)} \right) \forall n \geqslant 2
  \end{array} \label{licont}
\end{equation}
since
\begin{equation}
  \frac{\Gamma (b_n) \Gamma (b_{n - 1})}{\Gamma (a_{n + 1}) \Gamma (b_n + b_{n
  - 1} - a_{n + 1})} = \frac{\Gamma (2) \Gamma (2)}{\Gamma (1) \Gamma (2 + 2 -
  1)} = \frac{1}{2}
\end{equation}
and
\begin{equation}
  \frac{(b_n - a_{n + 1})_m (b_{n - 1} - a_{n - 1})_m}{(b_n + b_{n - 1} - a_{n
  + 1}) m!} = \frac{(1)_m (1)_m}{(2 + 2 - 1)_m m!} = \frac{2}{(m + 1) (m + 2)}
\end{equation}
The denominator parameters $\vec{2}^{+ 1 + m}_{n - 1}$ in (\ref{licont}) are
simply
\begin{equation}
  \begin{array}{lll}
    \vec{2}^{+ 1 + m}_{n - 1} & = \vec{2}_{n - 2}, 3 + m & = \underbrace{2,
    \ldots, 2}_{n - 2}, 3 + m
  \end{array}
\end{equation}
The numbers ($m + 1) (m + 2$) are known as the {\tmem{oblong numbers}},
{\cite[A002378]{oeis}}. By simply setting $t = 1$ in (\ref{licont}) we get the
continuation from $\zeta^F (n) \rightarrow \zeta^F (n + 1) \forall n \geqslant
1$
\begin{equation}
  \begin{array}{ll}
    \zeta_{}^F (n) =_{n + 1} F_n \left( \begin{array}{l}
      \vec{1}_{n + 1}\\
      \vec{2}_n
    \end{array}   \right) & = \sum_{m = 0}^{\infty} \left( \frac{_n F_{n - 1}
    \left( \begin{array}{l}
      \vec{1}_n\\
      \vec{2}^{+ m + 1}_{n - 1}
    \end{array} \right)}{(m + 1) (m + 2)} \right) \label{zetacont} \forall n
    \geqslant 2
  \end{array}
\end{equation}
The justification in saying that $\tmop{Li}_n^F (t)$ and $\zeta^F (n)$ are
continued to $\tmop{Li}_{n + 1}^F (t)$ and $\zeta^F (n + 1)$ comes from the
fact that the first term in the summand of the continuation (\ref{licont})
from $\tmop{Li}_{n - 1}^F (t) \rightarrow \tmop{Li}_n^F (t$) is contiguous to
\ $\tmop{Li}_{n - 1}^F (t)$, that is, $_n F_{n - 1} \left( \begin{array}{l}
  \vec{1}_n\\
  \vec{2}^{+ 1}_{n - 1}
\end{array} | t \right)$ is contiguous to $\tmop{Li}_{n - 1}^F (t) = \text{$_n
F_{n - 1} \left( \begin{array}{l}
  \vec{1}_n\\
  \vec{2}^{}_{n - 1}
\end{array} | t \right)$}$. The continuation formula (\ref{licont}) gives
interesting answers for $n = 0$ and $n = 1$ which suggest an alternative to
``the analytic continuation'' of $\zeta \left( t \right)$ which is different
from the usual $\frac{1}{1 - 2^{- t}} \sum_{n = 0}^{\infty} (2 n + 1)^{- t}$.
We have
\begin{equation}
  \begin{array}{llll}
    \zeta_{}^F (0) & = \sum_{m = 0}^{\infty} \left( \frac{_0 F_1 \left(
    \begin{array}{l}
      \\
      m + 3
    \end{array} \right)}{(m + 1) (m + 2)} \right) &  & \\
    & = \sum_{m = 0}^{\infty} \frac{- \left( I_{m + 1} \left( 2 \right) m +
    I_{m + 1} \left( 2 \right) - I_m \left( 2 \right) \right) \Gamma \left( m
    + 3 \right)}{(m + 1) (m + 2)}  &  & \\
    & = I_0 (2) - 1 &  & \\
    & \approx 1.2795853023360 &  & \\
    \zeta_{}^F (1) & = \sum_{m = 0}^{\infty} \left( \frac{_1 F_1 \left(
    \begin{array}{l}
      1\\
      m + 3
    \end{array} \right)}{(m + 1) (m + 2)} \right) &  & \\
    & = \sum_{m = 0}^{\infty} \frac{e \left( \Gamma \left( m + 3 \right) -
    \Gamma \left( m + 2, 1 \right) m - 2 \Gamma \left( m + 2, 1 \right)
    \right)}{(m + 1) (m + 2)}  &  & \\
    & = \tmop{Ei} (1) - \gamma &  & \\
    & \approx 1.3179021514544 &  & 
  \end{array}
\end{equation}
where $\tmop{Ei} (x)$ is the {\tmem{exponential integral}}
{\cite[6.9.2]{htf1}}
\begin{equation}
  \begin{array}{ll}
    \tmop{Ei} (x) & = \gamma - \frac{\ln (x^{- 1})}{2} + \frac{\ln (x)}{2} +
    \sum_{k = 1}^{\infty} \frac{x^k}{k \Gamma (k + 1)}\\
    & = \gamma - \frac{\ln (x^{- 1})}{2} + \frac{\ln (x)}{2} + x_2 F_2 \left(
    \begin{array}{ll}
      1 & 1\\
      2 & 2
    \end{array} | x \right)
  \end{array}
\end{equation}
and $\Gamma \left( a, z \right)$ is the {\tmem{incomplete Gamma function}}
\begin{equation}
  \Gamma \left( a, z \right) = \Gamma \left( z \right) - \frac{z^a _1 F_1
  \left( \begin{array}{l}
    a\\
    a + 1
  \end{array} | - z \right)}{a}
\end{equation}
So we have the ``hypergeometrically continued'' values $\zeta_{}^F (0) = I_0
(2) - 1$ and $\zeta_{}^F (1) = \tmop{Ei} (1) - \gamma$ whereas the ``real''
values are $\zeta \left( 0 \right) = - \frac{1}{2}$ and $\zeta \left( 1
\right) = \infty$. In terms of reciprocal probability we have
\begin{equation}
  \begin{array}{ll}
    \zeta^F \left( 0 \right)^{- 1} & \cong 78.15\%\\
    \zeta^F \left( 1 \right)^{- 1} & \cong 75.88\%
  \end{array}
\end{equation}

\subsubsection{$\tmop{Li}_1^F (t) \rightarrow \tmop{Li}_2^F (t)$ and
$\zeta_1^F (t) \rightarrow \zeta_2^F (t)$}

The continuation $\zeta_n^F (t)$ from $n = 1 \rightarrow 2$ via
(\ref{zetacont}) is straightforward
\begin{equation}
  \begin{array}{ll}
    \text{$\zeta^F (2)$} & =_3 F_2 \left( \begin{array}{lll}
      1 & 1 & 1\\
      2 & 2 & 
    \end{array}   \right)\\
    & = \sum_{m = 0}^{\infty} \left( \frac{_2 F_1 \left( \begin{array}{ll}
      1 & 1\\
      & 3 + m
    \end{array} \right)}{(m + 1) (m + 2)} \right)\\
    & = \sum_{m = 0}^{\infty} \left( \frac{\sum_{k = 0}^{\infty} \frac{\Gamma
    (m + 3) \Gamma (k + 1)}{\Gamma (m + k + 3)}}{(m + 1) (m + 2)} \right) \\
    & = \sum_{m = 0}^{\infty} \frac{1}{(m + 1) (m + 2)} \frac{(m + 2)}{(m +
    1)}\\
    & = \sum_{m = 0}^{\infty} \frac{1}{(m + 1)^2}\\
    & = \frac{\pi^2}{6}
  \end{array}
\end{equation}
The continuation of $\tmop{Li}^F_1 (t)$ to $\tmop{Li}^F_2 (t)$ is a bit more
complicated
\begin{equation}
  \begin{array}{ll}
    \tmop{Li}_2^F \left( t \right) & =_3 F_2 \left( \begin{array}{lll}
      1 & 1 & 1\\
      & 2 & 2
    \end{array} | t   \right)\\
    & = \sum_{m = 0}^{\infty} \left( \frac{_2 F_1 \left( \begin{array}{ll}
      1 & 1\\
      & m + 3
    \end{array} | t \right)}{(m + 1) (m + 2)} \right)\\
    & = \sum_{m = 0}^{\infty} r_2 (m, t)
  \end{array}
\end{equation}
then $r_2 (m, t)$ is given by

\begin{equation}
  \begin{array}{ll}
    r_{2 } (m, t) & = \frac{\sum_{n = 0}^m \frac{(- 1)^{n + 1} \Gamma (m + 2)
    (\Psi (m - n + 1) - \Psi (m + 2)) (- 1)^m e^{\psi (m + 2)} t^n}{\Gamma (n
    + 2) \Gamma (m - n + 1)}}{(m + 1) e^{\psi (m + 2)} t^{m + 1}}\\
    & - \frac{(t - 1)^{m + 1} t^{- 2 - m} \ln (1 - t)}{ m + 1}
  \end{array}
\end{equation}

so $\tmop{Li}_2^F (t)$ is equal to
\begin{equation}
  \text{$\tmop{Li}_2^F (t)$=} \sum_{m = 0}^{\infty} \frac{\sum_{n = 0}^m
  \frac{(- 1)^{n + 1} \Gamma (m + 2) (\Psi (m - n + 1) - \Psi (m + 2)) (- 1)^m
  e^{\psi (m + 2)} t^n}{\Gamma (n + 2) \Gamma (m - n + 1)}}{(m + 1) e^{\psi (m
  + 2)} t^{m + 1}} - \frac{(t - 1)^{m + 1} t^{- 2 - m} \ln (1 - t)}{m + 1}
\end{equation}
where $\psi (m) = \ln (\tmop{lcm} (1, 2, 3, \ldots, m))$ is Chebyshev's
function of the 2nd kind (\ref{cheb2}) and $\Psi (m$) is the digamma function
\begin{equation}
  \Psi (x) = \frac{\mathd}{\mathd x} \ln (\Gamma (x)) =
  \frac{\frac{\mathd}{\mathd x} \Gamma (x)}{\Gamma (x)} \label{digamma}
\end{equation}

\subsubsection{$\zeta^F (2) \rightarrow \zeta^F (3)$}

The continuation from $\zeta^F (2)$ to $\zeta^F (3)$ via (\ref{zetacont}) is
carried out like so
\begin{equation}
  \begin{array}{ll}
    \text{$\zeta^F (3)$} & =_4 F_3 \left( \begin{array}{llll}
      1 & 1 & 1 & 1\\
      & 2 & 2 & 2
    \end{array}   \right)\\
    & = \sum_{m = 0}^{\infty} \left( \frac{_3 F_2 \left( \begin{array}{lll}
      1 & 1 & 1\\
      & 2 & 3 + m
    \end{array} \right)}{(m + 1) (m + 2)} \right)\\
    & = \sum_{m = 0}^{\infty} r_3 (m)
  \end{array} \label{z3f}
\end{equation}
Each term in the summand $r_3 (m)$ has the form $\frac{\zeta (2)}{m + 1} + q_3
(m)$ where of course $\zeta (2) = \frac{\pi^2}{6}$ and $q_3 (m)$ is a rational
function of $m$ which follows a 3rd order linear recurrence
equation{\cite[8.2]{ab}} given by

\begin{equation}
  \begin{array}{ll}
    q_3 (m) & = q_3 \left( m + 1 \right) \left( m^3 + 8 \hspace{0.25em} m^2 +
    21 \hspace{0.25em} m + 18 \right)\\
    & + q_3 \left( m + 2 \right) \left( - 2 \hspace{0.25em} m^3 - 20
    \hspace{0.25em} m^2 - 67 \hspace{0.25em} m - 75 \right)\\
    & + q_3 \left( m + 3 \right) \left( m^3 + 12 \hspace{0.25em} m^2 + 48
    \hspace{0.25em} m + 64 \right)
  \end{array}
\end{equation}
\begin{equation}
  q_3 (m) = \left\{ \begin{array}{ll}
    - 1 & m = 0\\
    - \frac{5}{8} & m = 1\\
    - \frac{49}{108} & m = 2
  \end{array} \right.
\end{equation}
The solution to which is given by
\begin{equation}
  \begin{array}{ll}
    q_3 (m) & = \frac{\Psi^{(1)} (m + 2) - \zeta \left( 2 \right)}{m + 1}
  \end{array}
\end{equation}
so the summand $r_3 (m)$ is
\begin{equation}
  \begin{array}{ll}
    r_3 (m) & = \frac{\zeta (2)}{m + 1} + q_3 (m) = \frac{\Psi^{(1)} (2 +
    m)}{m + 1} \label{r3}
  \end{array}
\end{equation}
Thus (\ref{z3f}) is also equal to
\begin{equation}
  \begin{array}{ll}
    \zeta^F (3) & = \sum_{m = 0}^{\infty} \frac{\Psi^{(1)} (m + 2)}{m + 1}
  \end{array}
\end{equation}
Thus
\begin{equation}
  \begin{array}{ll}
    r_3 (m) & = \frac{\Psi^{(1)} (m + 2)}{m + 1}\\
    & = \frac{\zeta (2, m + 2)}{m + 1}\\
    & = \frac{\pi^2}{6} - \sum_{k = 1}^{m - 1} \frac{1}{k^2}\\
    & = \frac{\sum_{k = 1}^{\infty} \frac{1}{(k + m - 1)^2}}{m + 1}\\
    & = \frac{_3 F_2 \left( \begin{array}{lll}
      1 & m + 2 & m + 2\\
      & m + 3 & m + 3
    \end{array} \right)}{(m + 1) (m + 2)^2} \\
    & = \frac{_3 F_2 \left( \begin{array}{lll}
      1 & 1 & 1\\
      & 2 & 3 + m
    \end{array} \right)}{(m + 1) (m + 2)} 
  \end{array}
\end{equation}
The first 10 terms of $\{q_3 (m) : m = 0 \ldots 9\}$ are
\begin{equation}
  \left[ - 1, - \frac{5}{8}, - \frac{49}{108}, - \frac{205}{576}, -
  \frac{5269}{18000}, - \frac{5369}{21600}, - \frac{266681}{1234800}, -
  \frac{1077749}{5644800}, - \frac{9778141}{57153600}, -
  \frac{1968329}{12700800} \right] \label{z3terms}
\end{equation}
The denominator of (\ref{z3terms}) appears to be {\cite[A119936]{oeis}}, the
least common multiple of denominators of the rows of a certain triangle of
rationals and the numerators are {\cite[A007406]{oeis}}, the numerator of
$\sum_{k = 1}^n \frac{1}{k^2}$ from (\ref{phi1}) which, according to a theorem
Wolstenholme, $p$ divides $\tmop{numer} (q_3 (p - 1$)) where $p \in
\mathbbm{P}$ is prime. {\cite{wolstconv}} {\cite{wolstdiv}} {\cite{wolstgen}}

\subsubsection{$\zeta^F (3) \rightarrow \zeta^F (4)$}

The continuation from $\zeta^F (3)$ to $\zeta^F (4)$ via (\ref{zetacont}) is
given by
\begin{equation}
  \begin{array}{ll}
    \text{$\zeta^F (4)$} & =_5 F_4 \left( \begin{array}{lllll}
      1 & 1 & 1 & 1 & 1\\
      2 & 2 & 2 & 2 & 
    \end{array} \right)\\
    & = \sum_{m = 0}^{\infty} \left( \frac{_4 F_3 \left( \begin{array}{llll}
      1 & 1 & 1 & 1\\
      & 2 & 2 & 3 + m
    \end{array} \right)}{(m + 1) (m + 2)} \right) \\
    & = \sum_{m = 0}^{\infty} r_4 (m)
  \end{array} \label{z3f}
\end{equation}
The summand $r_4 (m)$ has the form
\begin{equation}
  \begin{array}{ll}
    r_4 (m) & = a (t, m) - b (t, m) - \frac{H (m + 1) \tmop{Li}_2 (t)}{(m + 1)
    t} + \frac{\tmop{Li}_3 (t)}{(m + 1) t}
  \end{array}
\end{equation}
where $a (t, m)$ is an $(m + 1)$-th degree polynomial and $b (t, m)$ is a $(m
+ 2)$-th degree polynomial(the determination of which is left to an excercise
for the reader or the topic of another article, but is readily obtained with
the help of Maple{\cite{maple15pg}}), and $H (n)$ is the $n$-th Harmonic
number
\begin{equation}
  \begin{array}{ll}
    H (n) & = \sum_{i = 1}^n \frac{1}{n}\\
    & = \Psi (n + 1) + \gamma\\
    & = \sum_{k = 1}^{\infty} \frac{n}{k^2 + k n}\\
    & = \frac{n}{n + 1} _3 F_2 \left( \begin{array}{lll}
      1 & 1 & n + 1\\
      & 2 & n + 2
    \end{array} \right)
  \end{array}
\end{equation}
The polynomial $b (t, m)$ vanishes when $t = 1$. An interesting set of
formulas for $\zeta (4)$ is
\begin{equation}
  \begin{array}{ll}
    \zeta (4) & = \sum_{n = 1}^{\infty} \frac{\Psi^{(2)} (n + 1) + 2 \zeta
    (3)}{2 n \left( n + 1 \right)}\\
    & = \sum_{n = 1}^{\infty} \frac{\Psi^{(2)} (n + 1) + 2 \sum_{m =
    0}^{\infty} \frac{\Psi^{(1)} (m + 2)}{m + 1}}{2 n \left( n + 1 \right)}\\
    & = \frac{\pi^4}{90}
  \end{array}
\end{equation}

\section{Appendix}

\subsection{$\tmop{The} \tmop{Generalized} \tmop{Hypergeometric}
\tmop{Function} :_p F_q$}

The Pochhammer symbol is defined according to
\begin{equation}
  \begin{array}{ll}
    (n)_k & = \frac{\Gamma (n + k)}{\Gamma (n)}
  \end{array} \label{pochhammer}
\end{equation}
The generalized hypergeometric function
{\cite{hypergeometricRepresentation}}{\cite[4.1]{whittaker-watson}} is defined
as an infinite sum of quotients of finite products of Pochhammer symbols

\begin{equation}
  \begin{array}{ll}
    \begin{array}{l}
      _p F_q \left( \begin{array}{l}
        a_1, \ldots, a_p\\
        b_1, \ldots, b_q
      \end{array} | t \right)
    \end{array} & = \sum_{k = 0}^{\infty} \frac{t^k}{k!} \frac{\prod_{i = 1}^p
    (a_i)_k}{\prod_{j = 1}^q (b_j)_k} \label{hypergeom}
  \end{array}
\end{equation}

The function $_p F_q$ is said to be $k$-balanced {\cite{hypercont}} if the sum
of the denominator parameters $b_1 \ldots b_p$ minus the sum of the numerator
parameters $a_1 \ldots a_{p + 1}$ is an integer.
\begin{equation}
  k = \tmop{bal} (_p F_q) = \sum_{n = 1}^q b_n - \sum_{n = 1}^p a_n
  \label{balance}
\end{equation}
The value $k$ is the characteristic exponent of the hypergeometric
differential equation at unit argument which is equal to the maximum root of
the corresponding indicial equation and so determines the behaviour of the
function near this point. A $1$-balanced function is said to be
\tmtextit{Saalsch\"utzian}. {\cite[2.1.1]{Fslater}}

\subsubsection{The Differential Equation and Convergence}

The function $_p F_q$ converges when
\begin{equation}
  \begin{array}{l}
    \left\{ \begin{array}{lll}
      p \leqslant q & \forall |t| \neq \infty & \\
      p = q + 1 & \forall |t| < 1 & \\
      \{p = q + 1 : \tmop{bal} (_p F_q) \geqslant 1\} & \forall |t| = 1 & \\
      p > q + 1 & \forall t = 0 & 
    \end{array} \right.
  \end{array} \label{Fconvergence}
\end{equation}
where $\tmop{bal} (_p F_q) = \sum_{n = 1}^q b_n - \sum_{n = 1}^p a_n$ is the
balance of the parameters (\ref{balance}). The differential equation solved by
$_p F_q \tmop{is} \tmop{of} \tmop{order}  \text{max(p,q+1)}$

\begin{equation}
  \left( \theta_t \prod_{j = 1}^q (\theta_t + b_j - 1) - t \prod_{i = 1}^p
  (\theta_t + a_i) \right) f (t) = 0 \label{hode}
\end{equation}

where $f (t) =_p F_q \left( \begin{array}{l}
  a_1, \ldots, a_p\\
  b_1, \ldots, b_q
\end{array} | t \right)$ and $\theta_t = t \frac{\mathd}{\mathd t}$ is the
differential operator. When $p = q + 1$ (\ref{hode}) has the form
\begin{equation}
  a_0 f (t) + t^q  \frac{\mathd}{\mathd t^{q + 1}} f (t) + \sum_{n = 1}^q t^{n
  - 1} (t a_n - b_n) \frac{\mathd}{\mathd t^n} f (t) = 0
\end{equation}
\cite[4.2]{whittaker-watson}{\cite[Ch3]{hypersum}}{\cite[44-46]{sf}}
{\cite[2.1.2]{Fslater}}

\subsubsection{Contiguous Functions and Linear Relations}\label{contig}

Any two hypergeometric functions $ {_p F_q} (a_{\ldots}, b_{\ldots} ; z)$ and
$_p F_q (c_{\ldots}, d_{\ldots} ; z)$ are said to be contiguous if all $p + q$
pairs of parameters $(a_1, c_1), \ldots, (a_p, c_p), (b_1, d_1), \ldots, (b_q,
d_q)$ are equal except for one pair which differs only by 1. There are $2 p +
q$ linearly independent relations between the $(2 p + 2 q)$ functions
contiguous to $_p F_q (a_{\ldots}, b_{\ldots} ; z)$ where the relations are
linear functions of $z$ and polynomial functions of the parameters
$a_{\ldots}, b_{\ldots}$. When any $\left\{ a_i = a_j : i \neq j\} \right.$ or
$\left\{ b_i = b_j : i \neq j\} \right.$ in $_p F_q (a_{\ldots}, b_{\ldots} ;
z)$ there will fewer unique contiguous functions than if all the parameters
were unique since the hypergeometric function is invariant with respect to the
ordering of parameters. {\cite[2.2.1]{Fslater}} {\cite[48]{sf}}
{\cite{hypergeometricRepresentation}} {\cite[4.3]{htf1}} {\cite{hyperlie}}
{\cite{pfq-contig}}

\subsection{Other Special Functions}

\subsubsection{Polygamma $\Psi^{(n)} (x)$}

The polygamma function is the $n$-th derivative of the digamma (\ref{digamma})
function
\begin{equation}
  \Psi^{(n)} (x) = \frac{\mathd}{\mathd x^n} \Psi (x)
\end{equation}
and is defined as an infinite sum, a Hurwitz zeta function (\ref{hurwitz}), or
a hypergeometric function when $x$ is positive integer
\begin{equation}
  \begin{array}{ll}
    \Psi^{(n)} (x) & = \left\{ \begin{array}{ll}
      \left( \sum_{k = 1}^{\infty} \frac{1}{k} - \frac{1}{k + x - 1} \right) -
      \matheuler & n = 0\\
      \sum_{k = 0}^{\infty} - \frac{n! (- 1)^n}{(k + x)^{n + 1}} & n \geqslant
      1
    \end{array} \right.\\
    & = \left\{ \begin{array}{ll}
      \left( \frac{x - 1}{x} _3 F_2 \left( \begin{array}{lll}
        1 & 1 & x\\
        & 2 & x + 1
      \end{array} \right) \right) - \matheuler & n = 0\\
      \frac{n! (- 1)^{n + 1}}{x^{n + 1}} _{n + 2} F_{n + 1} \left(
      \begin{array}{ll}
        1 & \vec{x}_{n + 1}\\
        \overrightarrow{(1 + x)}_{n + 1} & 
      \end{array} \right) & n \geqslant 1
    \end{array} \right.\\
    & = \left\{ \begin{array}{ll}
      \Psi (x) & n = 0\\
      (- 1)^{n + 1} n! \zeta (x, n + 1) & n \geqslant 1
    \end{array} \right.
  \end{array} \label{phiinf}
\end{equation}
or as a finite sum when $x$ is a positive integer and $n = 1$
{\cite[1.16]{htf1}}
\begin{equation}
  \begin{array}{ll}
    \Psi^{(1)} (x) & = \frac{\pi^2}{6} - \sum_{k = 1}^{x - 1} \frac{1}{k^2}
    \label{phi1}
  \end{array}
\end{equation}

\subsection{Notation}

\begin{equation}
  \begin{array}{ll}
    \mathbbm{Z} & \{\ldots, - 2, - 1, 0, 1, 2, \ldots\}\\
    \mathbbm{N} & \{1, 2, 3, \ldots .\}\\
    \mathbbm{N}_{1^-} & \{\ldots, - 3, - 2, - 1\}\\
    \mathbbm{N}_0 & \{0, 1, 2, 3, \ldots .\}\\
    \mathbbm{N}_{0^-} & \{\ldots, - 3, - 2, - 1, 0\}
  \end{array} \label{notation}
\end{equation}

\end{document}